\documentclass[12pt]{amsart}

\usepackage{amsmath,amsfonts,amsthm,amssymb}
\usepackage{hyperref}
\usepackage{tikz}
\usetikzlibrary{arrows.meta, positioning, calc}

\usepackage{xcolor}

\usepackage{a4wide}
\usepackage{enumerate}
\usepackage{multicol, multirow, longtable}
\newtheorem{theorem}{Theorem}[section]

\newtheorem{theo}{Theorem}

\newtheorem{coro}[theo]{Corollary}

\theoremstyle{remark}

\newtheorem*{claim*}{Claim}

\theoremstyle{definition}
\newtheorem{Definition}[theorem]{Definition}

\numberwithin{equation}{section}

\newcommand{\R}{\ensuremath{\mathbb{R}}}
\newcommand{\HH}{\ensuremath{\mathbb{H}}}
\newcommand{\g}[1]{\ensuremath{\mathfrak{#1}}}

\DeclareMathOperator{\Ric}{Ric}

\DeclareMathOperator{\rank}{rank}
\DeclareMathOperator{\Sym}{Sym}

\newcommand{\Cc}{\ensuremath{\mathcal{C}}}
\newcommand{\Dd}{\ensuremath{\mathcal{D}}}

\renewcommand{\sp}[1]{\ensuremath{\mathsf{Sp}_{#1}}}

\newcommand{\lieG}{{\mathfrak{g}}}
\newcommand{\lieH}{{\mathfrak{h}}}

\newcommand{\lieP}{{\mathfrak{p}}}

\usepackage{parskip}

\begin{document}
\title{Positive intermediate curvatures and Ricci flow}

\author[D.~Gonz\'alez-\'Alvaro]{David Gonz\'alez-\'Alvaro}
\address{David Gonz\'alez-\'Alvaro\\Universidad Polit\'ecnica de Madrid, Spain.}
\email{david.gonzalez.alvaro@upm.es}

\author[M.~Zarei]{Masoumeh Zarei}
\address{Masoumeh Zarei\\Mathematisches Institut, Universit\"at M\"unster, 
Einsteinstr. 62, 48149 M\"unster, Germany.}
\email{masoumeh.zarei@uni-muenster.de}

\begin{abstract}
We show that, for any $n\geq 2$, there exists a homogeneous space of dimension $d=8n-4$ with metrics of $\Ric_{\frac{d}{2}-5}>0$ if $n\neq 3$ and $\Ric_6>0$ if $n=3$ which evolve under the Ricci flow to metrics whose Ricci tensor is not $(d-4)$-positive. Consequently, Ricci flow does not preserve a range of curvature conditions that interpolate between positive sectional and positive scalar curvature. This extends a theorem of B\"ohm and Wilking in the case of $n=2$.
\end{abstract}


\maketitle

\section{Introduction}

Given a Riemannian manifold $(M,g)$, it is a fundamental problem to understand how the metric $g$ and its curvature properties evolve under the Ricci flow $\frac{\partial}{\partial t}g(t)=-2\Ric(g(t))$ for $g(0)=g$, see e.g. \cite{Ni14,BK19,Ba21}. In this note we study the evolution of metrics satisfying curvature conditions which interpolate between positive sectional and positive scalar curvature.

\begin{theo}\label{T:Main_New}
For any $n\geq 2$ there exists a homogeneous space of dimension $d=8n-4$ with metrics of $\Ric_{\frac{d}{2}-5}>0$ if $n\neq 3$ and $\Ric_6>0$ if $n=3$ which evolve under the Ricci flow to metrics whose Ricci tensor is not $(d-4)$-positive. 
\end{theo}

Recall that the conditions $\Ric_{k}>0$ interpolate between positive sectional, when $k=1$, and positive Ricci curvature, when $k=d-1$, where $d$ denotes the dimension of the manifold. Furthermore, $k$-positivity of the Ricci tensor, i.e. the sum of its  $k$ smallest eigenvalues being positive, are conditions which interpolate between positive Ricci, when $k=1$, and positive scalar curvature, when $k=d$. An immediate consequence of Theorem~\ref{T:Main_New} is that wide  ranges of curvature conditions are not preserved under the Ricci flow.  In order to be more precise, we define the following cones of algebraic curvature operators for each $d$. 

Let $\Sym^{2}_{b}(\wedge^{2}\mathbb{R}^{d})$ denote the space of algebraic curvature operators, which are defined as the elements in $\Sym^{2}(\wedge^{2}\mathbb{R}^{d})$ satisfying the first Bianchi identity. Recall that for $R\in \Sym_b^{2}(\wedge^{2}\mathbb{R}^{d})$, the sectional curvature of a $2$-plane spanned by orthornormal vectors $x,y$ is defined as $\langle R(x\wedge y), x\wedge y\rangle$. More generally, the $\Ric_k^{R}$ curvature of a flag spanned by orthonormal vectors $\{x,y_1,\dots,y_k\}$ with $k\leq d-1$ is defined as $\sum_i\langle R(x\wedge y_i), x\wedge y_i\rangle$. Let $\Ric^R$ be the Ricci tensor associated to $R\in \Sym_b^{2}(\wedge^{2}\mathbb{R}^{d})$. We consider the following cones of algebraic curvature operators:
\begin{align*}
\Cc_k &= \{ R\in \Sym_b^{2}(\wedge^{2}\mathbb{R}^{d})\; :\; \Ric^{R}_k>0 \text{ for all flags} \}, \\
\Dd_k &=\{ R\in \Sym_b^{2}(\wedge^{2}\mathbb{R}^{d})\; :\; \text{the sum of the $k$ smallest eigenvalues of $\Ric^R$ is positive}\}.
\end{align*}
In this terminology we have the chain of inclusions
$$
\Cc_1\subset\dots\subset\Cc_k\subset\Cc_{k+1}\subset\dots\subset\Cc_{d-1}=\Dd_1\subset\dots\subset\Dd_k\subset\Dd_{k+1}\subset\dots\subset\Dd_d.
$$
Now, we say that a Riemannian manifold $M$ satisfies the curvature condition $C_k$ (resp. $D_k$) if, under an identification $T_pM\cong \R^d$ via an orthonormal basis of $T_pM$, the Riemannian curvature operator $R_p$ satisfies $R_p\in\Cc_k$ (resp. $R_p\in\Dd_k$) for all $p\in M$. Note that $C_1$, $C_{d-1}=D_1$ and $D_d$ correspond to positive sectional, Ricci and scalar curvatures, respectively.

\begin{coro}\label{coro}
\begin{enumerate}
\item There is no curvature condition between $C_{\frac{d}{2}-5}$ and $D_{d-4}$ which is invariant under the Ricci flow on closed simply connected manifolds of dimension $d$, for any $d=8n-4\geq 12$ and $d\neq 20$. 

\item  There is no curvature condition between $C_{6}$ and $D_{16}$ which is invariant under the Ricci flow on  closed simply connected manifolds of dimension $20$. 
\end{enumerate}
\end{coro}

In order to put Theorem~\ref{T:Main_New} and Corollary~\ref{coro} into perspective, let us recall that by  work of Hamilton \cite{Ha82}, positive scalar curvature is preserved under the Ricci flow in every dimension. Moreover, both positive sectional and positive Ricci curvature are preserved in dimension $3$. It is then natural to ask whether any other curvature conditions are preserved in higher dimensions. Now we review some results that provide partial answers to this question.

There exist examples in infinitely many dimensions where positive Ricci curvature is not preserved \cite{AN16}, as well as examples in dimensions $4,6,12,24$ where positive sectional curvature is not preserved \cite{BK23,BW07,CW15}. In this context, the example in dimension $12$ is specially relevant to us. Namely, B\"ohm and Wilking \cite{BW07} show that there are metrics of positive sectional curvature on the Wallach space $W^{12}=\sp{3}/\sp{1}^3$ for which positive Ricci curvature is not preserved. Indeed, their proof shows that such metrics evolve to metrics whose Ricci tensor is not $8$-positive.

The Wallach space $W^{12}$ can be seen as the first member in the following infinite family of spaces. For $n\geq 2$ we denote by $P_n$ the $(8n-4)$-dimensional homogeneous space
$$
P_n = \sp{n+1}/(\sp{n-1}\times\sp{1}\times\sp{1}).
$$
Observe that $P_n$ is the projectivized tangent bundle $\mathbb{P}_\HH T\mathbb{HP}^n$ of $\mathbb{HP}^n$. In \cite[Theorem~G]{DGM23} it is proven that $P_n$ admits an $\sp{n+1}$-invariant metric of $\Ric_{k(n)}>0$ for the following values of $k(n)$:
$$
k(n)=
\begin{cases}
4n-7, & n\neq 3,\\
6, & n=3.\\
\end{cases}
$$
We can now state a refined version of Theorem~\ref{T:Main_New}, which extends the result of B\"ohm and Wilking for $P_2$ to all of the spaces $P_n$.

\begin{theorem}\label{T:Main}
For any $n\geq 2$, there exist $\sp{n+1}$-invariant metrics of $\Ric_{k(n)}>0$ on $P_n$ that evolve under the Ricci flow to metrics whose Ricci tensor has $8n-8$ negative eigenvalues.
\end{theorem}

Note that by taking Riemannian products, one can extend Theorem~\ref{T:Main_New} and Corollary~\ref{coro} to other dimensions, albeit at the expense of covering less curvature conditions relative to the dimension. For example, products $P_n\times\mathbb{S}^m$ with round spheres of dimension $m\geq 2$ carry metrics of $\Ric_k>0$ with $k=\max\{ k(n)+m, 8n-3 \}$ which evolve under the Ricci flow to metrics whose Ricci tensor is not $(8n-8)$-positive. Note that $k>\frac{d}{2}$, where $d=8n-4+m$ denotes the dimension of $P_n\times\mathbb{S}^m$. The smallest $k$'s are reached when $m$ is roughly half of the dimension of $P_n$, in which case $k$ converges to $\frac{2d}{3}$ as $n$ tends to infinity. Similar observations can be made for other products like $P_n\times P_m$ or $P_n\times T^m$, where $T^m$ denotes the $m$-torus.

Let us emphasize that for $d$-dimensional manifolds, the curvature bounds $\Ric_k>0$ with $k\leq\frac{d}{2}$ seem to be significantly stronger than those with $k>\frac{d}{2}$. As noted in the case of $P_n\times\mathbb{S}^m$, the fairly trivial construction of taking Riemannian products produces manifolds of $\Ric_k>0$ only for $k>\frac{d}{2}$. Recently, Reiser and Wraith have refined in \cite{RW23, RW25} surgery and bundle-like techniques to construct a vast number of manifolds of $\Ric_k>0$ for various $k>\frac{d}{2}$. In particular, they have shown that Gromov's Betti number bound fails in the case of $\Ric_{[\frac{d}{2}]+2}>0$ in any dimension $d\geq 5$.  Finally, we mention that the bound $k\leq\frac{d}{2}$ also appears in several structural results of manifolds of $\Ric_k>0$, see \cite{GW22,KM24,Mo22,Xi97}. We refer to \cite{Mo} for a collection of publications and preprints concerning the curvature conditions $\Ric_k>0$.

Unfortunately, there is no fixed notation nor terminology in the literature for the curvature conditions considered  in this note. 
For example, $k$-positivity of the Ricci tensor is denoted by $\Ric_k>0$ in the article \cite{CW22}. We refer to \cite[Section~2.2]{DGM23} or \cite[p.~5]{CW22} for further information. In order to avoid any confusion, we include the precise definitions and notation that we use.

\begin{Definition}
Let $M$ be a $d$-dimensional Riemannian manifold and let $k\leq d-1$. We say that $M$ is of \emph{positive $k^{th}$-intermediate Ricci curvature}, to be denoted by $\Ric_k>0$, if for every point $p\in M$ and every set of $k+1$ orthonormal tangent vectors $v,e_1,\dots,e_k\in T_p M$, the sum of the sectional curvatures of the planes spanned by $v,e_i$ is positive, i.e.
\[
\sum_{i=1}^k \sec(v,e_i)>0.
\]
\end{Definition}

\begin{Definition}
Let $M$ be a $d$-dimensional Riemannian manifold and let $k\leq d$. We say that the Ricci tensor of $M$ is \emph{$k$-positive} if the sum of the $k$ smallest eigenvalues of the Ricci tensor is positive at all points.
\end{Definition}

It is easy to see that $\Ric_{k}>0$ implies $\Ric_{k+1}>0$, and it follows from the definition that $k$-positivity of the Ricci tensor implies $(k+1)$-positivity. In the context of Theorem~\ref{T:Main_New} and Theorem~\ref{T:Main}, the implications between the relevant curvature conditions on a manifold of dimension $d=8n-4$ are:
\begin{align*}
\sec>0\Rightarrow\Ric_{k(n)}>0\Rightarrow\Ric_{8n-5} >0\Leftrightarrow\Ric>0\Rightarrow \{\text{Ricci tensor is $k$-positive for } k\geq 1\}.
\end{align*}

The relations between the corresponding cones of curvature operators are:
$$
\Cc_1\subset\Cc_{k(n)}\subset\Cc_{8n-5}\subset\Dd_k,\qquad k\geq 1.
$$
Note that with this notation, Theorem~\ref{T:Main} in fact implies that Ricci flow evolves certain metrics with curvature operator in $\Cc_{k(n)}$ to metrics with curvature operator not in $\Dd_{8n-8}$.

\subsection*{Acknowledgements} We would like to thank Anusha M. Krishnan for helpful conversations. We are grateful to Christoph B\"ohm and Philipp Reiser for useful comments on a preliminary version of this manuscript. This work was done while the first author was visiting the University of M\"unster. The first author wishes to thank the University of M\"unster for providing excellent working conditions. We thank an anonymous referee for suggesting improvements to the exposition.

The first author was supported by grants PID2021-124195NB-C31 and PID2021-124195NB-C32 from the Agencia Estatal de Investigación and the Ministerio de Ciencia e Innovación (Spain). The second author was supported  by Deutsche Forschungsgemeinschaft (DFG, German Research Foundation) under Germany's Excellence Strategy EXC 2044-390685587, Mathematics M\"unster: Dynamics-Geometry-Structure, and by  DFG grant ZA976/1-1 within the Priority Program SPP2026 ``Geometry at Infinity''.

\section{Proof of Theorem~\ref{T:Main}}

In this section we first recall various known results about the spaces $P_n$ and then we prove Theorem~\ref{T:Main}.

\subsection{Preliminaries}\label{SS:preliminaries}

The homogeneous spaces $P_n$ from the introduction belong to a class of homogeneous spaces called \emph{generalized Wallach spaces}. These spaces are characterized as compact homogeneous spaces $G/H$ whose isotropy representation decomposes into a direct sum $\lieP=\lieP_1\oplus\lieP_2\oplus\lieP_3$ of three $Ad(H)$-invariant irreducible modules satisfying $[\lieP_i,\lieP_i] \subset \lieH$  for $i=1, 2, 3$, where $\g{h}$ is the Lie algebra of $H$. Generalized Wallach spaces have been recently classified \cite{Niko16,Niko21}. 

From now on we assume that $G/H$ is a generalized Wallach space and that the modules $\g{p}_i$ are pairwise inequivalent. Denote by $Q$ the bi-invariant inner product on $\lieG$ given by the negative of the Killing form of $G$. Any $G$-invariant metric on $G/H$ is determined by an $Ad(H)$-invariant inner product $g$ on $\g{p}$ of the form
$$g=x_1Q|_{\lieP_1}+x_2Q|_{\lieP_2}+x_3Q|_{\lieP_3},$$ 
for some $x_1, x_2, x_3>0$. 
The Ricci tensor of any such metric is given by
$$
\Ric(g)=r_1x_1Q|_{\lieP_1}+r_2x_2Q|_{\lieP_2}+r_3x_3Q|_{\lieP_3},
$$
for certain numbers $r_1,r_2,r_3$, which can be computed by the formula (see \cite[Section~3.1]{St22})
\begin{equation}\label{eq:ri}
r_i=\frac{1}{2x_i} +\frac{a_i}{2}\left(\frac{x_i}{x_jx_k}-\frac{x_k}{x_ix_j}-\frac{x_j}{x_ix_k} \right).
\end{equation}
Here $a_1,a_2,a_3$ are rational numbers depending on $G/H$ and can be found in \cite[Table~1]{Niko16}. By choosing a $g$-orthonormal basis adapted to the decomposition $\lieP=\lieP_1\oplus\lieP_2\oplus\lieP_3$ we see that the eigenvalues of the Ricci tensor are $r_1,r_2,r_3$ with  multiplicities $d_1,d_2,d_3$, respectively, where $d_i=\dim\g{p}_i$. 

Note that for the spaces $P_n$, since $\rank G=\rank H$, the irreducible modules $\g{p}_i$ are pairwise inequivalent. Therefore, the discussion above is applicable to $P_n$. Moreover, we have
$$
d_1=d_2=4(n-1),\qquad d_3=4,
$$
and the constants $a_i$ are equal to
$$
a_1=a_2=\frac{1}{2(n+2)},\qquad a_3=\frac{n-1}{2(n+2)}.
$$

\subsection{The metrics of $\Ric_{k(n)}>0$}\label{SS:metrics}

Consider the nested inclusions
$$
\sp{n-1}\times\sp{1}\times\sp{1}<\sp{n-1}\times\sp{2}< \sp{n+1}.
$$
This leads to a homogeneous fibration
$$
\sp{2}/(\sp{1}\times\sp{1}) \to P_n \to \sp{n+1}/(\sp{n-1}\times\sp{2}),
$$
where the base space is a quaternionic $2$-plane Grassmannian and the fiber is the sphere $\mathbb{S}^4$. A basic method to construct metrics on $P_n$ is to start with a normal homogeneous metric (i.e. a metric with $x_1=x_2=x_3$) and rescale it in the direction of the fiber. Such metrics will be called \emph{submersion metrics}. Since the tangent space to the fiber $\mathbb{S}^4$ corresponds to $\lieP_3$, it follows that submersion metrics correspond to triples $(x_1,x_2,x_3)$ with $x_1=x_2\neq x_3$. In \cite[Section~5]{DGM23} it is proven that any metric on $P_n$ with $x_1=x_2>x_3$ is of $\Ric_{k(n)}>0$.

\subsection{Normalized Ricci flow}

The sign of each eigenvalue of the Ricci tensor is preserved under rescalings of the metric. Thus, for our purposes, it is enough to look at the \emph{normalized Ricci flow} equation
$$
\frac{\partial}{\partial t}g(t)=-2\Ric(g(t))+ 2g(t)\frac{S_{g(t)}}{n},
$$
where $S_{g(t)}$ denotes the scalar curvature of $g(t)$. The normalized Ricci flow $g(t)$ of an initial metric $g=g(0)$ preserves its isometry group and hence in our case we have:
\begin{equation}\label{eq:metrics_flow}
g(t)=x_1(t)Q|_{\lieP_1}+x_2(t)Q|_{\lieP_2}+x_3(t)Q|_{\lieP_3},
\end{equation}
for some functions $x_1(t), x_2(t), x_3(t)>0$. Then the normalized Ricci flow reduces to a system of ordinary differential equations
$$
x_1'=f(x_1,x_2,x_3),\quad x_2'=g(x_1,x_2,x_3),\quad x_3'=h(x_1,x_2,x_3).
$$
In the case of generalized Wallach spaces $G/H$ the expressions for the functions $f,g,h$ can be found in \cite[page~27]{AANS14}. 

The upshot of using the normalized Ricci flow is that the volume stays constant. Hence, we may assume that $x_1^{d_1}x_2^{d_2}x_3^{d_3}=1$ along the flow, which allows to express $x_3$ in terms of $x_1$ and $x_2$. In the case of $P_n$ we get $x_3=(x_1x_2)^{-(n-1)}$. Thus we are left with a system of just two equations:
$$
x_1'=f(x_1,x_2,(x_1x_2)^{-(n-1)}),\quad x_2'=g(x_1,x_2,(x_1x_2)^{-(n-1)}).
$$
Altogether, by using the concrete expressions for $f$ and $g$ in the case of $P_n$, we get the system
\begin{equation}
\begin{aligned}
x_1' &=-1-\frac{x_1}{2(n+2)}\left(x_1^{n}x_2^{n-2} - x_1^{n-2}x_2^{n} - \frac{1}{x_1^n x_2^n} \right) + x_1 B \\
x_2' &=-1-\frac{x_2}{2(n+2)}\left( - x_1^{n}x_2^{n-2} + x_1^{n-2}x_2^{n} - \frac{1}{x_1^n x_2^n} \right) + x_2 B,
\end{aligned}
\label{eq:ODEsystem_x}
\end{equation}
where
$$
B=\left(  2(n+2)\left(\frac{1}{x_1} +  \frac{1}{x_2} + \frac{x_1^{n-1}x_2^{n-1}}{n-1}\right) - x_1^{n}x_2^{n-2} - x_1^{n-2}x_2^{n} - \frac{1}{x_1^n x_2^n} \right)\frac{n-1}{2(n+2)(2n-1)}.
$$

\subsection{Proof of Theorem~\ref{T:Main}}

We follow the strategy of B\"ohm-Wilking in \cite{BW07} to prove Theorem~\ref{T:Main}. First we make the following change of coordinates:
$$
\varphi = x_1 + x_2, \qquad \psi=x_1 - x_2.
$$
Note that $\varphi>0$ and $\varphi>\psi$.

From \eqref{eq:ri} and \eqref{eq:metrics_flow} we can write $r_i(t)$ of the metric $g(t)$ in the new coordinates as follows:
\begin{equation}\label{eq:r1}
\begin{aligned}
r_1 &=\frac{1}{\varphi + \psi}+\frac{\varphi\psi(\varphi^2-\psi^2)^{n-2}}{4^{n-1}(n+2)} - \frac{4^{n-1}}{(n+2)(\varphi^2-\psi^2)^{n}},\\
&\\
r_2&=\frac{1}{\varphi - \psi}+\frac{\varphi\psi(\varphi^2-\psi^2)^{n-2}}{4^{n-1}(n+2)} - \frac{4^{n-1}}{(n+2)(\varphi^2-\psi^2)^{n}}, \\
&\\
r_3&=\frac{(\varphi^2-\psi^2)^{n-1}}{2\cdot 4^{n-1}} - \frac{(n-1)(\varphi^2-\psi^2)^{n-2}}{4^{n-1}(n+2)} + \frac{4^{n-1}(n-1)}{(n+2)(\varphi^2-\psi^2)^{n}}.
\end{aligned}
\end{equation}

Furthermore, System~\eqref{eq:ODEsystem_x} becomes:
\begin{equation}
\begin{aligned}
\varphi' =& -2  +\frac{ (\varphi^2 - \psi^2)^{n-2}}{4^{n-1}(n+2)(2n-1)}\left(3\varphi^3 - (6n-1)\varphi\psi^2 \right)\\
& +\frac{4^nn\varphi}{2(n+2)(2n-1)(\varphi^2 - \psi^2)^{n}} + \frac{n-1}{2n-1}\left(\frac{4\varphi^2}{\varphi^2 - \psi^2}\right)  \\
&\\
\psi' =& \psi\left(  \frac{(\varphi^2 - \psi^2)^{n-2}}{4^{n-1}(n+2)(2n-1)}\left((-4n+5)\varphi^2 -(2n+1)\psi^2 \right)\right.  \\
 & \left. +\frac{4^nn}{2(n+2)(2n-1)(\varphi^2 - \psi^2)^{n}} + \frac{n-1}{2n-1}\left(\frac{4\varphi}{\varphi^2 - \psi^2}\right)  \right).
\end{aligned}
\label{eq:ODEsystem_new_coordinates}
\end{equation}

For each initial condition corresponding to a homogeneous metric there is some $T\in(0,\infty]$ such that System~\eqref{eq:ODEsystem_new_coordinates} has a solution for all $t\in [0,T]$ (see \cite[Theorem~3.1]{Sb22} for more information on the existing time $T$). 

The first observation is that if $(\varphi,\psi)$ is a solution with $\psi(t_0)=0$ for some $t_0\in [0,T]$, then $\psi(t)=0$ for all $t\in [0,T]$. In other words, the Ricci flow of a submersion metric on $P_n$ stays within the class of submersion metrics. To show it, take $\psi(t)\equiv 0$, so that System~\eqref{eq:ODEsystem_new_coordinates} reduces to the equation
\begin{equation}\label{eq:psi0}
\varphi' = \frac{1}{2n-1}\left(-2+\frac{3\varphi^{2n-1}}{4^{n-1}(n+2)}+\frac{4^n n}{2(n+2)\varphi^{2n-1}} \right).
\end{equation}
This equation has a solution $\varphi$ for each initial condition. Thus $(\varphi,0)$ is a solution of System~\eqref{eq:ODEsystem_new_coordinates}. By uniqueness we are done.

An immediate implication of this observation is that if $\psi(t_0)<0$ for some $t_0\in [0,T]$, then $\psi(t)<0$ for all $t\in [0,T]$. 

When $\psi\equiv 0$ and $\varphi(0)$ is large enough, it follows from \eqref{eq:psi0} that $\varphi'(t)>1$ for all $t\in [0,T]$. Under these conditions it is easy to see that Ricci flow preserves $r_i(t)>0$ for $i=1,2,3$.

In order to find metrics for which the Ricci flow does not preserve $r_i>0$ for some $i$, we consider initial metrics which are perturbations of submersion metrics. More precisely, we assume $\varphi(0)=N$ is large and $\psi(0)<0$ with $|\psi(0)|$ very small. These metrics are very close to submersion metrics and moreover $x_3=(x_1x_2)^{-(n-1)}<x_1$; thus, since $\Ric_k>0$ is an open condition for any $k$, such metrics will be of $\Ric_{k(n)}>0$ as discussed in Section~\ref{SS:metrics}.

Now we fix $\varphi(0)=N$ sufficiently large and $\psi(0)<0$ with $|\psi(0)|$ very small (to be fixed later). These conditions imply, by analyzing System~\eqref{eq:ODEsystem_new_coordinates}, that $\varphi'(t)>1$ and $\psi'(t)>0$ for all $t\in [0,T]$. Because $0<\frac{1}{\varphi'(t)}<1$, the equation $h'(t)=\frac{1}{\varphi'(h(t))}$ with initial condition $h(0)=0$ has a solution $h$ defined for all $t\in [0,\infty)$. We use this map to do the following reparametrization:
$$
\tilde\varphi(t):=\varphi(h(t)),\qquad \tilde\psi(t):=\psi(h(t)).
$$
This yields a new system for $\tilde\varphi '$ and $\tilde\psi '$, stated below, whose solutions are reparametrizations of solutions of System~\eqref{eq:ODEsystem_new_coordinates}. One advantage of this reparametrization is that, even if some solutions of \eqref{eq:ODEsystem_new_coordinates} only exist up to a finite time, the corresponding solutions to the new system  are defined for all $t\in [0,\infty)$.

By abuse of notation we simply write $\varphi$, $\psi$ instead of $\tilde\varphi$, $\tilde\psi$; as only the behavior of these maps for large $t$ is relevant to us. System~\eqref{eq:ODEsystem_new_coordinates} now reads as:
\begin{align*}
\varphi' =& 1 \\
\psi' =& \psi\frac{\left(  \frac{(\psi^2 - \varphi^2)^{n-2}}{4^{n-1}(n+2)(2n-1)}\left((-4n+5)\varphi^2 -(2n+1)\psi^2 \right)  +\frac{4^nn}{2(n+2)(2n-1)(\varphi^2 - \psi^2)^{n}} + \frac{n-1}{2n-1}\left(\frac{4\varphi}{\varphi^2 - \psi^2}\right)  \right)}{-2  +\frac{ (\varphi^2 - \psi^2)^{n-2}}{4^{n-1}(n+2)(2n-1)}\left(3\varphi^3 - (6n-1)\varphi\psi^2 \right)+\frac{4^nn\varphi}{2(n+2)(2n-1)(\varphi^2 - \psi^2)^{n}} + \frac{n-1}{2n-1}\left(\frac{4\varphi^2}{\varphi^2 - \psi^2}\right)}.
\end{align*}
Clearly we have $\varphi(t)=t+N$. Next we shall obtain an upper bound for $\psi'$ for large $t$. Observe:
$$
\lim_{t\to\infty} \frac{\psi'\varphi}{\psi} =\frac{-4n+5}{3}.
$$
Thus, for each $\eta > \frac{4n-5}{3}$ there is a time $t_0(\eta)$ such that 
$$
\psi '\leq -\eta \psi\frac{1}{\varphi}, \qquad\text{for all } t\geq t_0(\eta).
$$
We take $\eta=\frac{4n}{3}-1$ and denote the corresponding $t_0(\eta)$ simply by $t_0$. By integrating both sides we get that there is some constant $C>0$ such that
$$
\psi \leq  -C\frac{1}{\varphi^\eta}, \qquad\text{for all } t\geq t_0.
$$
Now we fix $\psi(0)<0$ with $|\psi(0)|$ sufficiently small so that at time $t=0$ we have the reverse inequality $\psi(0) \geq  -C\frac{1}{\varphi(0)^\eta}$. Thus, we may assume that at time $t_0$ we have an equality $\psi(t_0)=  -C\frac{1}{\varphi(t_0)^\eta}$. Hence we get the inequality
$$
\psi(t) \leq  \psi(t_0)\varphi(t_0)^\eta\frac{1}{\varphi(t)^\eta}, \qquad\text{for all } t\geq t_0.
$$
By multiplying both sides by $\varphi(t)^{2n-2}$ and noting that $2n-2 -\eta=\frac{2n}{3}-1>0$, we derive that 
\begin{equation}\label{eq:limit}
\lim_{t\to\infty} \psi(t)\varphi(t)^{2n-2} \leq \lim_{t\to\infty} \psi(t_0)\varphi(t_0)^\eta\varphi(t)^{2n-2-\eta} = -\infty.
\end{equation}
Using the formula \eqref{eq:r1} for $r_1$ we can write
$$
r_1 \varphi = \frac{\varphi}{\varphi + \psi}+\frac{\varphi^{2n-2}\psi}{4^{n-1}(n+2)}\left(\frac{\varphi^2-\psi^2}{\varphi^2}\right)^{n-2} - \frac{4^{n-1}\varphi}{(n+2)(\varphi^2-\psi^2)^{n}}.
$$
Using \eqref{eq:limit} it clearly follows that $\lim_{t\to\infty} r_1 \varphi=-\infty$. Because $\varphi$ is always positive it follows that $r_1(t)$ is negative for all $t$ sufficiently large. This shows that, for the corresponding solutions of System~\eqref{eq:ODEsystem_new_coordinates}, $r_1(t)<0$ for $t$ sufficiently close to $T$. The same arguments show that $r_2(t)<0$ for $t$ sufficiently close to $T$. Recall that $r_1(t)$ and $r_2(t)$ are eigenvalues of the Ricci tensor of $g(t)$ with the same multiplicity $d_1=d_2=4(n-1)$. This proves Theorem~\ref{T:Main}.


\end{document}